\documentclass[12pt, journal, draftclsnofoot, onecolumn]{IEEEtran}
\ifCLASSOPTIONcompsoc
\else
\fi

\ifCLASSINFOpdf
\else
\fi
\usepackage[parfill]{parskip}    
\usepackage{amsmath,amssymb,latexsym,enumerate, mathrsfs}
\usepackage{graphicx}
\usepackage{graphicx}
\usepackage{array}
\usepackage{times}
\usepackage{hyperref}
\newtheorem{lemma}{\bf Lemma}

\newtheorem{definition}{Definition}
\newtheorem{proposition}{Proposition}

\def\qed{\hfill $\Box$}
\def\etal{\mbox{et al.}}

\begin{document}
%
\title{A brief remark on the topological entropy for linear switched systems}

\author{Getachew~K.~Befekadu~\IEEEmembership{}
\IEEEcompsocitemizethanks{\IEEEcompsocthanksitem G. K. Befekadu is with the Department
of Electrical Engineering, University of Notre Dame, Notre Dame, IN 46556, USA.\protect\\
E-mail: gbefekadu1@nd.edu\protect\\
Version - February 25, 2013.}
\thanks{}}

\markboth{}%
{Shell \MakeLowercase{\textit{et al.}}: Bare Advanced Demo of IEEEtran.cls for Journals}
\IEEEcompsoctitleabstractindextext{%
\begin{abstract}
In this brief note, we investigate the topological entropy for linear switched systems. Specifically, we use the Levi-Malcev decomposition of Lie-algebra to establish a connection between the basic properties of the topological entropy and the stability of switched linear systems. For such systems, we show that the topological entropy for the evolution operator corresponding to a semi-simple subalgebra is always bounded from above by the negative of the largest real part of the eigenvalue that corresponds to the evolution operator of a maximal solvable ideal part.
\end{abstract}

}

\maketitle

\IEEEdisplaynotcompsoctitleabstractindextext

%
\IEEEpeerreviewmaketitle

\section{Introduction} \label{sec:Intro}
In the past, the notions of measure-theoretic entropy and topological entropy have been intensively studied in the context of measure-preserving transformations or continuous maps (e.g., see references \cite{Wal82}, \cite{Sin94} and \cite{Dow11} for the review of entropy in ergodic theory as well as in dynamical systems). For instance, Adler \etal\, (in the paper \cite{AdlKM65}) introduced the notion of topological entropy as an invariant conjugacy, which is an analogue to the notion of measure-theoretic entropy, for measuring the rate at which a continuous map in a compact topological space generates initial-state information. Subsequently, Bowen and Dinaburg, in the papers \cite{Din70} and \cite{Bow71} respectively, gave a weak, but equivalent, definition of topological entropy for continuous maps that later led to proofs for connecting this notion of entropy with that of measure-theoretic entropy (e.g., see also \cite{Goodm71} or \cite{Goodm72} for additional discussions). 

With the emergence of networked control systems (e.g., see \cite{AntBa07}), these notions of entropy have found renewed interest in the research community (e.g., see \cite{NaiEMM04}, \cite{Sav06} and \cite{ColKa09}). Notably, Nair \etal\, \cite{NaiEMM04} have introduced the notion of topological feedback entropy, which is based on the ideas of \cite{AdlKM65}, to quantify the minimum rate at which deterministic discrete-time dynamical systems generate information relevant to the control objective of set-invariance. More recently, the notion of (controlled)-invariance entropy has been studied for continuous-time control systems in \cite{ColKa09} and \cite{ColKa11} based on the metric-space technique of \cite{Bow71}. It is noted that such an invariance entropy provides a measure of the smallest growth rate for the number of open-loop control functions that are needed to confine the states within an arbitrarily small distance from a given compact subset of the system state space.

On the other hand, several results have been established to characterize the stability and/or the performance of switched systems using Lie-brackets -- where feasible, one can consider the Lie-algebra generated by the constituent systems (a matrix Lie-algebra in the linear case or a Lie-algebra of vector fields in the general nonlinear case) and use this Lie-algebra to verify or establish the stability of switched systems (e.g., see \cite{LibHM99}, \cite{AgrL01}, \cite{Lib03} and references therein for a review of switched systems). Note that if the constituent systems are linear and stable, it should be noted that the switched system remains stable under an arbitrary switching if the Lie-algebra is nilpotent (e.g., see \cite{Gur95}) or a compact semi-simple subalgebra (e.g., see \cite{AgrL01}). Moreover, it should be noted that each of these classes of Lie-algebras strictly contains the others (e.g., see \cite{Che62}, \cite{Lib03}). 

\section{Preliminaries}\label{sec:prelim}
\subsection{Switched systems}
Consider a Lie-algebra (over a real $\mathbb{R}$) that is generated by the matrices $A_p \in \mathbb{R}^{n \times n}$, $p \in \mathcal{P}\equiv  \{1, 2, \ldots, N\}$ and further identified with
\begin{align}
\mathfrak{g} = \bigl\{A_p \colon p \in \mathcal{P} \bigr\}_{LA}. \label{Eq2.1}
\end{align}
Let $\mathfrak{g}= \mathfrak{m} \oplus \mathfrak{h}$ be the Levi-Malcev decomposition of the Lie-algebra, where $\mathfrak{m}$ is the radical (i.e., the maximal solvable ideal part) and $ \mathfrak{h}$ is the semi-simple subalgebra part. Then, we can rewrite the matrices $A_p$ as
\begin{align}
A_p = A_p^{\mathfrak{m}} + A_p^{\mathfrak{h}}, \label{Eq2.2}
\end{align}
with $A_p^{\mathfrak{m}} \in \mathfrak{m}$ and $A_p^{\mathfrak{h}} \in \mathfrak{h}$ for $p \in \mathcal{P}$.

Next, we consider the following family of systems
\begin{align}
\dot{x}(t) = A_{\sigma(t)} x(t), \quad x(0)=x_0, \label{Eq2.3}
\end{align}
where $\sigma \colon [0,\, \infty) \to \mathcal{P}$ is a piecewise constant switching function.\footnote{In this brief note, switching that is infinitely fast (i.e., chattering) is not considered. Moreover, we assume that the matrices $A_p \in \mathbb{R}^{n \times n}$, $\forall p \in \mathcal{P}$, are stable.}

Let $\Phi(t,\,0)$ or simply $\Phi(t)$ (assuming that the initial time $t_0$ is zero) be the {\em evolution operator} for the family of systems in \eqref{Eq2.3} and observe that
\begin{align}
\dot{\Phi}(t) &= A_{\sigma(t)} \Phi(t), \notag \\
                     &= \bigl(A_{\sigma(t)}^{\mathfrak{m}} + A_{\sigma(t)}^{\mathfrak{h}}\bigr) \Phi(t), \label{Eq2.4}
\end{align}
with $\sigma \colon [0,\, \infty) \to \mathcal{P}$. 

Then, we state the following well-known result that will be useful in the sequel.
\begin{lemma} [{\bf Levi-Malcev Decomposition}] \label{lem:L1} 
The evolution operator $\Phi(t)$ can be decomposed as follow 
\begin{align}
\Phi(t) = \Phi^{\mathfrak{h}}(t) \Phi^{\mathfrak{m}}(t), \label{Eq2.5}
\end{align}
where 
\begin{align}
\dot{\Phi}^{\mathfrak{h}}(t) = A_{\sigma(t)}^{\mathfrak{h}} \Phi^{\mathfrak{h}}(t), \quad  \Phi^{\mathfrak{h}}(0) = I, \label{Eq2.6}
\end{align}
and
\begin{align}
\dot{\Phi}^{\mathfrak{m}}(t) = \biggm(\bigl(\Phi^{\mathfrak{h}}(t)\bigr)^{-1} A_{\sigma(t)}^{\mathfrak{m}} \Phi^{\mathfrak{h}}(t) \biggm) \Phi^{\mathfrak{m}}(t), \quad  \Phi^{\mathfrak{m}}(0) = I.  \label{Eq2.7}
\end{align}
\end{lemma}

The proof follows the same lines of argument as that of \cite{Che62}.

{\em Proof:} Note that if we differentiate equation~\eqref{Eq2.5}, i.e., the evolution operator $\Phi(t)$, with respect to time and also make use of the relations in \eqref{Eq2.6} and \eqref{Eq2.7}, then we have
\begin{align}
\frac{d}{dt}\bigm(\Phi^{\mathfrak{h}}(t) \Phi^{\mathfrak{m}}(t)\bigm) &= \frac{d}{dt}\bigm(\Phi^{\mathfrak{h}}(t)\bigm) \Phi^{\mathfrak{m}}(t) + \Phi^{\mathfrak{h}}(t) \frac{d}{dt}\bigm(\Phi^{\mathfrak{m}}(t)\bigm), \notag \\
&= A_{\sigma(t)}^{\mathfrak{h}} \Phi^{\mathfrak{h}}(t) \Phi^{\mathfrak{m}}(t) + \Phi^{\mathfrak{h}}(t)  \frac{d}{dt}\bigm(\Phi^{\mathfrak{m}}(t)\bigm), \notag \\
&= A_{\sigma(t)}^{\mathfrak{h}} \Phi^{\mathfrak{h}}(t) \Phi^{\mathfrak{m}}(t) + \Phi^{\mathfrak{h}}(t) \biggm(\bigm(\Phi^{\mathfrak{h}}(t)\bigm)^{-1} A_{\sigma(t)}^{\mathfrak{m}} \Phi^{\mathfrak{h}}(t) \biggm) \Phi^{\mathfrak{m}}(t), \notag \\
&= \bigm(A_{\sigma(t)}^{\mathfrak{h}} + A_{\sigma(t)}^{\mathfrak{m}}\bigm) \Phi^{\mathfrak{h}}(t) \Phi^{\mathfrak{m}}(t), \notag \\
&= A_{\sigma(t)} \Phi^{\mathfrak{h}}(t) \Phi^{\mathfrak{m}}(t). \label{Eq2.8}
\end{align}
\qed

Next, we introduce the following notion of stability for the family of systems in \eqref{Eq2.3}.

\begin{definition}
The switched system in \eqref{Eq2.3} is said to be globally uniformly exponentially stable (GUES), if there exist positive numbers $M$ and $\lambda$ such that the solutions of \eqref{Eq2.3} satisfy
\begin{align}
\vert x(t)\vert \le M \exp(-\lambda t)\vert x(0)\vert, \quad \forall t \ge 0. \label{Eq2.9}
\end{align}
\end{definition}

\subsection{Topological entropy for switched systems}
We start by providing the definition of topological entropy for switched linear systems that corresponds to the semi-simple subalgebra part (e.g., see \cite{Bow71} or \cite{Wal82} for additional discussions on the topological entropy for continuous transformations).
\begin{definition}
A set $\mathscr{F}$ is {\em $(T,\,\epsilon)$-spanning} another set $\mathscr{K}$ (with respect to $\Phi^{\mathfrak{h}}(t)\equiv\Phi_t^{\mathfrak{h}}$) if there exists $y \in \mathscr{F}$ for each $x \in \mathscr{K}$ such that 
\begin{align}
 \sup_{\substack{t \in [0,\, T]}} \bigl\Vert \Phi_t^{\mathfrak{h}}x - \Phi_t^{\mathfrak{h}}y\bigr\Vert \le \epsilon,  \label{Eq2.10}
\end{align}
where $\epsilon$ is a positive real number. 
\end{definition}
For a compact subset $\mathscr{K} \subset \mathscr{X}$ (where $\mathscr{X}$ is a compact $n$-dimensional $C^{\infty}$ manifold), let $r(T, \epsilon, \mathscr{K}, \Phi_t^{\mathfrak{h}})$ be the smallest cardinality of any subset $\mathscr{F} \subset \mathscr{X}$ that $(T,\,\epsilon)$-spans the set $\mathscr{K}$.\footnote{Note that the compactness of $\mathscr{X}$ implies that there exit finite $(T,\,\epsilon)$-spanning sets.} Then, we have the following properties for $r(T, \epsilon, \mathscr{K}, \Phi_t^{\mathfrak{h}})$.
\begin{enumerate} [(i)]
\item Clearly $r(T, \epsilon, \mathscr{K}, \Phi_t^{\mathfrak{h}}) \in [0,\,\infty)$.
\item If $\epsilon_1 < \epsilon_2$, then $r(T, \epsilon_1, \mathscr{K}, \Phi_t^{\mathfrak{h}}) \ge r(T, \epsilon_2, \mathscr{K}, \Phi_t^{\mathfrak{h}})$.
\end{enumerate}

\begin{definition}
The topological entropy for the switched linear system that corresponds to the semi-simple subalgebra is given by 
\begin{align}
h(\mathscr{K}, \Phi_t^{\mathfrak{h}}) =\,\lim_{\substack{\epsilon \searrow 0}} \biggm\{\limsup_{\substack{T \to \infty}} \frac{1}{T} \log r(T, \epsilon, \mathscr{K}, \Phi_t^{\mathfrak{h}})\biggm\}. \label{Eq2.11}
\end{align}
\end{definition}
Then, we have the following additional properties for $h(\mathscr{K}, \Phi_t^{\mathfrak{h}})$.
\begin{enumerate} [(i)]
\item $h(\mathscr{K}, \Phi_t^{\mathfrak{h}}) \in [0,\,\infty) \cup \{ \infty \}$.
\item If $\mathscr{K} = \bigcup_{l \in \{1,2,\ldots, L\}} \mathscr{K}_l$ with compact $\mathscr{K}_l$, then $h(\mathscr{K}, \Phi_t^{\mathfrak{h}}) = \max_{\substack{l \in \{1,2,\ldots, L\}}}h(\mathscr{K}_l, \Phi_t^{\mathfrak{h}})$.
\end{enumerate}

\section{Main result} \label{sec:main}
In the following, using the Levi-Malcev decomposition of the Lie-algebra, we establish a connection between the topological entropy for the evolution operator (that corresponds to the semi-simple subalgebra part of the Lie-algebra) and the stability of switched linear systems.

\begin{proposition} \label{prop:P1} 
Let $\mathscr{K} \subset \mathscr{X}$ be a compact subset. Suppose that the Levi-Malcev decomposition of the Lie-algebra corresponding to the matrices $A_p$ with $p \in \mathcal{P}$ is given by \eqref{Eq2.5}. Furthermore, if the topological entropy for the switched linear system corresponding to the semi-simple subalgebra part satisfies the following
\begin{align}
h(\mathscr{K}, \Phi_t^{\mathfrak{h}})  < - \bar{\lambda}_p^{\mathfrak{m}}, \quad \forall p \in \mathcal{P}, \label{Eq3.1}
\end{align}
where
\begin{align}
\bar{\lambda}_p^{\mathfrak{m}} = \max \biggm \{\operatorname{Re}\{\lambda\} \colon \lambda \in \operatorname{Sp}\bigl(A_p^{\mathfrak{m}}\bigr) \biggm\} , \quad p \in \mathcal{P}. \label{Eq3.2}
\end{align}
Then, the family of systems in \eqref{Eq2.3} are GUES.\footnote{$\operatorname{Sp}(A_p^{\mathfrak{m}})$ denotes the spectrum for the matrix $A_p^{\mathfrak{m}} \in \mathbb{R}^{n \times n}$.}
\end{proposition}
{\em Proof:} Note that if the evolution operator $\Phi(t)$ in \eqref{Eq2.2} admits a decomposition of the form in \eqref{Eq2.5}, then any $\mathfrak{m}$-valued solution $x_{\mathfrak{m}}(t)$ corresponding to $\Phi^{\mathfrak{m}}(t)$ satisfies the following
\begin{align}
\vert x_{\mathfrak{m}}(t)\vert \le \exp(\bar{\lambda}_p^{\mathfrak{m}} t)\vert x_{\mathfrak{m}}(0)\vert, \quad \forall t \in (0,\,T], \quad \forall p \in \mathcal{P}, \label{Eq3.3} 
\end{align}
with $\bar{\lambda}_p^{\mathfrak{m}}=\max \bigl \{\operatorname{Re}\{\lambda\} \colon \lambda \in \operatorname{Sp}\bigl(A_p^{\mathfrak{m}}\bigr)\bigr\}$ for $p \in \mathcal{P}$.

On the other hand, the characteristic Lyapunov exponent ${\lambda_p^{\mathfrak{h}}}^*$ for the evolution operator $\Phi_t^{\mathfrak{h}}$ is given by
\begin{align}
 {\lambda_p^{\mathfrak{h}}}^* =\, \lim_{\substack{t \to \infty}} \sup \frac{1}{t} \log\bigm\vert \operatorname{det}\bigl(\Phi^{\mathfrak{h}}(t)\bigr)\bigm\vert, \quad p \in \mathcal{P}, \label{Eq3.4}
\end{align}
where such information provides a lower bound for the topological entropy $h(\mathscr{K}, \Phi_t^{\mathfrak{h}})$ that corresponds to the semi-simple subalgebra part (e.g., see \cite{Yom87}, \cite{Pes77} or \cite{Kat80} for details on the relationships between Lyapunov exponents and entropy).\footnote{We remark that the topological entropy of a measure-preserving transformation always majorizes the measure-theoretic entropy with respect to any of its invariant probability measures (see also \cite{Wal82}).}

Then, the set of solutions for the family of systems in \eqref{Eq2.3} is exponentially bounded (i.e., the family of systems in \eqref{Eq2.3} are GUES), if the following condition holds (e.g., see also \cite{AgrBL12})
\begin{align}
\bar{\lambda}_p^{\mathfrak{m}} + {\lambda_p^{\mathfrak{h}}}^* < 0, \quad \forall p \in \mathcal{P}. \label{Eq3.5}  
\end{align}
Hence, this further implies the following
\begin{align}
\bar{\lambda}_p^{\mathfrak{m}} + h(\mathscr{K}, \Phi_t^{\mathfrak{h}}) < 0, \quad \forall p \in \mathcal{P}, \notag 
\end{align}
which completes the proof.
\qed

\end{document}